\documentclass{amsart}[12pt]
\def\doctype{preprint}

\usepackage{latexsym,amssymb}
\usepackage{fancyhdr}
\usepackage{lineno}
\usepackage{hyperref}

\newcommand{\cB}{\mathcal{B}}

\newcommand\Z{\mathbb{Z}}
\newcommand\F{\mathbb{F}}

\newcommand\cL{\mathcal{L}}

\newcommand{\comment}[1]{}

\numberwithin{equation}{section}

\fancypagestyle{titlepage}{
\fancyhead[R]{\doctype}
\fancyhead[CL]{}
\cfoot{\vspace{5pt} \thepage}
}

\newcommand\Section[1]{\section{\bf #1}}

\newtheoremstyle{theorem}
  {12pt}		  
  {0pt}  
  {\sl}  
  {\parindent}     
  {\bf}  
  {. }    
  { }    
  {}     
\theoremstyle{theorem}
\newtheorem{thm}{Theorem}[section]  
\newtheorem{lemma}[thm]{Lemma}     
\newtheorem{cor}[thm]{Corollary}
\newtheorem{prop}[thm]{Proposition}
\newtheorem{cons}[thm]{Construction}

\theoremstyle{definition}

\renewcommand\proof{\noindent\textsl{Proof. }}

\makeatletter
\renewcommand*\@maketitle{%
  \normalfont\normalsize
  \@adminfootnotes
  \@mkboth{\@nx\shortauthors}{\@nx\shorttitle}%
  \global\topskip42\p@\relax 
  \@settitle
  \ifx\@empty\authors \else {\vskip 1em
\vtop{\centering\shortauthors\@@par}} \fi
  \ifx\@empty\@date \else {\vskip 1em \vtop{\centering\@date\@@par}}\fi 
  \ifx\@empty\@dedicatory
  \else
    \baselineskip18\p@
    \vtop{\centering{\footnotesize\itshape\@dedicatory\@@par}%
      \global\dimen@i\prevdepth}\prevdepth\dimen@i
  \fi
  \@setabstract
  \normalsize
  \if@titlepage
    \newpage
  \else
    \dimen@34\p@ \advance\dimen@-\baselineskip
    \vskip\dimen@\relax
  \fi
} 
\renewcommand*\@adminfootnotes{%
  \let\@makefnmark\relax  \let\@thefnmark\relax
  \ifx\@empty\@subjclass\else \@footnotetext{\@setsubjclass}\fi
  \ifx\@empty\@keywords\else \@footnotetext{\@setkeywords}\fi
  \ifx\@empty\thankses\else \@footnotetext{%
    \def\par{\let\par\@par}\@setthanks}%
  \fi
\thispagestyle{titlepage}
}
\makeatother

\begin{document}

\title[PBD Dimension]{\large Pairwise balanced designs with prescribed \\ minimum dimension}

\author{Peter J.~Dukes}
\address{\rm Peter J.~ Dukes: Department of Mathematics and Statistics, 
University of Victoria, Victoria, BC, Canada}
\email{dukes@uvic.ca}

\author{Alan C.H.~Ling}
\address{\rm Alan C.H.~Ling: Department of Computer Science,
University of Vermont, Burlington, Vermont, U.S.A. 
} 
\email{aling@cems.uvm.edu}

\thanks{Research of Peter J.~Dukes is supported by NSERC}

\date{December 4, 2013}

\begin{abstract}
The dimension of a linear space is the maximum positive integer $d$ such that any $d$ of its points generate a proper subspace.  For a set $K$ of integers at least two, recall that a pairwise balanced design PBD$(v,K)$ is a linear space on $v$ points whose lines (or blocks) have sizes belonging to $K$.  We show that, for any prescribed set of sizes $K$ and lower bound $d$ on the dimension, there exists a PBD$(v,K)$ of dimension at least $d$ for all sufficiently large and numerically admissible $v$.
\end{abstract}

\maketitle
\hrule

\Section{Introduction}

An {\em incidence structure} is a triple $(X,\cL,\iota)$, where $X$ is a 
set of {\em points}, $\cL$ is a set of {\em lines}, and $\iota \subset 
X \times \cL$ is a set of {\em flags}.  We say $x \in X$ is {\em 
incident with} or simply \emph{on} $L \in \cL$ (and vice-versa) if and only if $(x,L) \in \iota$.  

A {\em linear space} is an incidence structure $(X,\cL,\iota)$ with the
property that every line is on at least two points and any two distinct
points are both on exactly one line.  In what follows $X$ (and 
hence $\cL$) are assumed finite. 
The trivial case in which all points are on the same line is not excluded by our definition but it is effectively ruled out.  (We would like to count lines as `subspaces' later on; apart from this, nontriviality can be assumed.)

Linear spaces appear in another context as {\em pairwise balanced 
designs} (or PBDs).  Specifically, if $v$ is a positive integer and 
$K \subset \Z_{\ge 2}:= \{2,3,4,\dots\}$ is a set of {\em block sizes}, a PBD$(v,K)$ 
consists of a $v$-set $X$, together with a set $\cB$ of {\em blocks}, where 
\begin{itemize}
\item
for each $B \in \cB$, we have $B \subset X$ with $|B| \in K$; and
\item
any two distinct elements of $X$ appear together in exactly one block.
\end{itemize}
Note that there are numerical constraints on $v$ given $K$.  First, the number of
pairs of distinct points must be expressible as a (nonnegative) integral linear combination of the number
of distinct pairs arising from blocks with sizes in $K$.  This leads to the \emph{global condition}
$$v(v-1) \equiv 0 \pmod{\beta(K)},\eqno{{\rm (global)}}$$
where $\beta(K):=\gcd\{k(k-1): k \in K\}$.  Also, deleting any point $x \in X$ from its incident blocks must partition
the remaining points.  That is, $v-1$ is an integral combination of $k-1$, $k \in K$.  This is the \emph{local condition}, namely
$$v-1 \equiv 0 \pmod{\alpha(K)},\eqno{{\rm (local)}}$$
where $\alpha(K):=\gcd\{k-1: k \in K\}$.

We can interchangeably discuss linear spaces and PBDs, identifying 
lines with the respective subsets of incident points as blocks.  So in what follows, notation such as $(X,\cB)$
is used for PBDs and the associated linear spaces; the incidence relation $\iota$ is seldom used from now on. 
However, we occasionally retain some terminology from 
linear spaces (i.e. points, lines, spaces) when discussing PBDs. 

Despite the similarity in the definitions, there is usually a difference in focus between the study of PBDs and linear spaces.
The former is usually approached with a fixed $K$ in mind, asking for which $v$ we have existence. 
The latter often concerns additional structures such as configurations or localizations at points.
Some features, such as parallelism, appear in both contexts.

Define $(X',\cB')$ as a  {\em subspace} (or  {\em subdesign}) of $(X,\cB)$ if $X' \subseteq X$ and $\cB' \subseteq \cB$.
That is, two distinct points in $X'$ must be covered by a unique block in $\cB'$.
Recalling that we permit trivial spaces, the points on a single line (block) always form a subspace.
As usual, a subspace on $X'$ is called {\em proper} if $X' \neq X$.

In $(X,\cB)$, the subspace {\em generated by} $Y \subset X$ is the unique 
minimal subspace containing $Y$.  
Let's denote this (set of points) by $\langle Y \rangle_{\cB}$, where
the subscript can be deleted if context is clear. 
On one hand, $\langle Y \rangle$ is the intersection of all subspaces containing $Y$. Alternatively, 
$\langle Y \rangle$ can be computed algorithmically starting from $Y$ by repeatedly including points on lines defined by existing points.

The {\em dimension} of a linear space 
is the maximum integer $d$ such that any set of $d$ points generates a 
proper subspace.  For instance, the subspace generated by any two points 
is the line containing them.  So every nontrivial linear space has dimension at least two.  
See \cite{D} and Chapter 7 of \cite{BB} for nice surveys of dimension in linear spaces.

It is unfortunate that the property of dimension has seldom made its way into 
the language of PBDs, and into design theory in general.  However, there 
are important exceptions to this.

Recall that a {\em Steiner triple system} is a PBD$(v,\{3\})$.  
It is well-known that 
Steiner triple systems on $v$ points exist if and only if $v \equiv 1$ or 
$3 \pmod{6}$.  
A {\em Steiner space} is defined to be a Steiner triple system of dimension at least 3.  
Teirlinck in \cite{T} nearly completely settled the existence of Steiner spaces. The result is 
that for $v \equiv 1$ or $3 \pmod{6}$ and $v \not\in \{51,67,69,145\}$,
there exists a Steiner space on $v$ points if and only if $v=15,27,31,39$, 
or $v \ge 45$.  The four undecided cases are still open, to the best of our knowledge.

Another important family of linear spaces, especially in design theory, is that of the affine spaces.
Let $q$ be a
prime power and $\mathbb{F}_q$ the finite field of order $q$.  Consider
the vector space $X=\F_q^{d}$ as points, together with all possible translates of subspaces 
$x+W \subseteq X$ as `flats'.
This forms the {\em affine space} AG$_d(q)$. 

Let $\cB$ be the set of all lines in AG$_d(q)$.  From basic linear algebra, we see that 
$(X,\cB)$ is a linear space (PBD) of dimension $d$, since $d$-point-generated subspaces correspond to proper flats (and some $d+1$ points generate the whole space).
There are $v=q^d$ points and every line has exactly $k=q$ 
points.  In other words, this is a PBD$(q^d,\{q\})$ of dimension $d$.

We hope for at least a small revival in the study of dimension in design
theory.  To this end, we present an existence theory that treats arbitrary block size(s) for all sufficiently large and admissible $v=|X|$.

\begin{thm}[Main theorem, full version]
\label{main}
Given $K \subseteq \Z_{\ge 2}$ and $d \in \Z_+$, there exists a PBD$(v,K)$ of dimension at least $d$ for all sufficiently large
$v$ satisfying {\rm (global)} and {\rm (local)}.
\end{thm}

When $K=\{k\}$, we have $\alpha(K)=k-1$ and $\beta(K)=k(k-1)$. Since this is often the case of primary interest, and for clarity of presentation, 
we first prove this case separately.  There is not much loss in economy, since most of the proof can be re-used to establish Theorem~\ref{main}.

\begin{thm}[Main theorem, weak version]
\label{main-k}
For $k \in  \Z_{\ge 2}$ and $d \in \Z_+$, there exists a PBD$(v,\{k\})$ of dimension at least $d$ for all sufficiently large
$v$ satisfying
\begin{eqnarray*}
v-1 &\equiv& 0 \pmod{k-1}; \text{and}\\
v(v-1) &\equiv& 0 \pmod{k(k-1)}.
\end{eqnarray*}
\end{thm}

Very broadly, the proofs proceed by applying some standard design-theoretic constructions to the affine space of dimension $d$, ensuring that the dimension stays preserved.  We introduce the needed background for the techniques in the next section.  Then, in Section 3, we carry out a certain sequence of constructions and prove that `intervals' of constructible values of $v$ eventually overlap.  Finally, we finish with a discussion of some related items, including a surprising connection with a problem in extremal graph theory.

\Section{Constructions}

First off, we state Wilson's famous `asymptotic' existence result for PBDs.

\begin{thm}[Wilson's Theorem, \cite{W}]
\label{Wilson}
Given $K \subseteq  \Z_{\ge 2}$, there exists $v_0$ such that a PBD$(v,K)$ exists for all $v \ge v_0$ satisfying {\rm (global)} and {\rm (local)}.
\end{thm}

Note that our main result simply says that one can demand a minimum dimension in this theorem.

The \emph{replication number} of a PBD$(v,\{k\})$ is the common number $r=\frac{v-1}{k-1}$ of blocks incident with each point.  The local necessary condition for $K=\{k\}$ amounts to $r  \equiv 0 \pmod{1}$.  The global condition is easily seen as equivalent to $r(r-1) \equiv 0 \pmod{k}$.  So we can restate Wilson's Theorem in terms of replication numbers.  This proves convenient in what follows.

\begin{cor}
\label{repl}
Given $k \ge 2$, there exists $r_0(k)$ such that a PBD with blocksize $k$ and replication number $r$ exists for all $r \ge r_0$ satisfying $r(r-1) \equiv 0 \pmod{k}$.
\end{cor}

It is clear that blocks of a PBD can be replaced by other PBDs.  That is, the existence of a PBD$(v,K)$ and, for each $k \in K$, a PBD$(k,L)$ implies the existence of a PBD$(v,L)$.  This construction, which is usually known as `breaking up blocks', respects dimension in a certain sense.

\begin{cons}
\label{bub}
Suppose there exists a PBD$(v,K)$ of dimension $d$ and, for each $k \in K$, any PBD$(k,L)$.  Then there exists a PBD$(v,L)$ of dimension $\ge d$.
\end{cons}

\proof
In the PBD$(v,K)$, say $(X,\cB)$, replace each block $B$ of size $k$ with a PBD$(k,L)$ on the points of $B$.  The result is a PBD$(v,L)$, say $(X,\cB_1)$.  It remains to check the dimension.  Suppose a set $Y$ of $d$ points is given.  They generate a proper subspace $X'$ in $(X,\cB)$ by hypothesis.  But this remains a subspace in $\cB_1$ after replacement of blocks by PBDs.  
\qed

A \emph{group divisible design} is a triple $(X,\Pi,\cB)$, where $X$ is a set of \emph{points}, $\Pi$ is a partition of $X$ into \emph{groups}, and $\cB$ is a set of \emph{blocks} such that 
\begin{itemize}
\item
a group and a block intersect in at most one point; and
\item
every pair of points from distinct groups is together in exactly one block.
\end{itemize}
We refer to this as a GDD or $K$-GDD, the latter emphasizing that the blocks have sizes in $K \subseteq \Z_{\ge 2}$. 
The \emph{type} of a GDD is the list of its group sizes.  When this list contains, say, $u$ copies of the integer $g$, this is abbreviated with `exponential notation' as $g^u$.  Another standard abbreviation is the use of $k$-GDD instead of $\{k\}$-GDD.

A \emph{transversal design} TD$(k,n)$ is a $k$-GDD of type $n^k$.  In this case, every block meets every group in one point.

We make use of two more important asymptotic existence results for the preceding objects.

\begin{thm}[Chowla, Erd\H{o}s, Strauss, \cite{CES}]
\label{asymptotic-td}
Given $k \ge 2$, there exists $n_0(k)$ such that a TD$(k,n)$ exists for all $n \ge n_0(k)$.
\end{thm}

\begin{thm}[Liu, \cite{Liu}]
\label{asymptotic-gdd}
Given $K \subseteq \Z_{\ge 2}$ and $g \in \Z_{+}$, there exists $u_0$ such that a $K$-GDD of type $g^u$ exists for all $u \ge u_0$ satisfying
\begin{eqnarray*}
g(u-1) &\equiv& 0 \pmod{\alpha(K)},\\
g^2 u(u-1) &\equiv& 0 \pmod{\beta(K)},
\end{eqnarray*}
where $\alpha$ and $\beta$ are as defined earlier.
\end{thm}

\emph{Remark}. One proof of Theorem~\ref{asymptotic-gdd} follows from edge-colored graph decompositions; see \S 8 of \cite{LW}.

We can regard a GDD $(X,\Pi,\cB)$ as the linear space (or PBD) $(X,\cB \cup \Pi)$, where groups and blocks are taken together to form lines.  Alternatively, we could consider the linear space $(X, \cB \cup \Pi_2)$, where $\Pi_2$ denotes the set of all pairs of distinct points from common groups.  Each of these interpretations allows one to talk about dimension for GDDs, the former being stronger (lower dimension) in general.
For our purposes, though, we prefer to take an even stronger notion for dimension.

Given a GDD, say $(X,\Pi,\cB)$, let's call a subspace $X' \subset X$ \emph{strong} if it intersects each group of $\Pi$ in either all points or no points.  A strong subspace is then proper if it is disjoint from at least one group.  (In practice, many groups will be missed.)  Correspondingly, the \emph{strong dimension} of a GDD is the maximum number of points which always generates a proper strong subspace.  With a PBD$(v,K)$ regarded as a $K$-GDD of type $1^v$, strong dimension coincides with ordinary dimension in this case.  On the other hand, the strong dimension of a transversal design is just 1, since two points from different groups generate a block of the TD, which in turn intersects all groups.

Given a PBD, say $(X,\cB)$, if we delete a point $x$ and all incident blocks $\cB_x$, the result is a GDD $(X \setminus \{x\}, \Pi_x, \cB \setminus \cB_x)$.  
Here, the group partition $\Pi_x$ is given by the (now missing) punctured lines $B \setminus \{x\}$, where $B \in \cB_x$.  Reversing this process, if we are given a GDD, say $(X,\Pi,\cB)$, we can add a point $\infty$ and replace groups with new blocks, all incident with $\infty$.  One might abbreviate this PBD by $(X^*,\cB^*)$, where $X^* = X \cup \{\infty\}$ and $\cB^* = \cB \cup \{X_i \cup \{\infty\} : X_i \in \Pi \}$.

\begin{lemma}
\label{add-pt}
If $(X,\Pi,\cB)$ has strong dimension $d$, then $(X^*,\cB^*)$ has dimension $\ge d$.
\end{lemma}

\proof
Consider a set $Y$ of $d$ points in $(X^*,\cB^*)$.  By hypothesis, $Y \setminus \{\infty\}$ is contained in a proper strong subspace $X'$ of $(X,\Pi,\cB)$.  After $\infty$ is included $X' \cup \{\infty\}$ becomes a proper subspace in $\cB^*$, since $\cB^*$ contains blocks $X_i \cup \{\infty\}$ for any group $X_i \subset X'$.
\qed

More generally, one can place PBDs, instead of single new blocks, on each extended group.  This is similar to Construction~\ref{bub}.

\begin{cons}
Suppose there exists a $K$-GDD on $v$ points with group sizes in $G$.  If, for each $g \in G$, there exists a PBD$(g+1,K)$, then there exists a PBD$(v+1,K)$.  
Furthermore, if the GDD has strong dimension $d$, then the resultant PBD has dimension $\ge d$.
\end{cons}

Rather than deleting a point, one could instead \emph{truncate} $x \in X$, replacing all blocks $B \in \cB_x$ by new blocks $B \setminus \{x\}$.  (New blocks of size 1 can be ignored.)  If the original space is a PBD or GDD, then so is the truncation.  It is a common design-theoretic technique to truncate several points from the same group of a GDD.  In this case, the modified blocks are only reduced in size by one.

\begin{cons}
\label{trunc}
If some, but not all, points of some group are truncated from a GDD, then its strong dimension does not decrease.
\end{cons}

\proof
Take a GDD $(X,\Pi,\cB)$ of strong dimension $d$, and truncate $Z \subset X$ from a common group.  Consider a set $Y$ of $d$ points in $X \setminus Z$.  By assumption, $Y$ is contained in a proper strong subspace $X' \subset X$.  Since no group has been deleted by the truncation, $X'$ remains a proper strong subspace in $X \setminus Z$.
\qed

Next is a powerful composition construction which played a key role in the proof of Theorem~\ref{Wilson}.

\begin{cons}[Wilson's fundamental construction]
\label{wfc}
Suppose there exists a `master' GDD $(X,\Pi,\cB)$, where $\Pi=\{X_1,\dots,X_u\}$.  Let $\omega:X \rightarrow \{0,1,2,\dots\}$, assigning nonnegative weights to each point in such a way that for every $B \in \cB$ there exists an `ingredient' $K$-GDD of type $\omega(B):=[\omega(x) \mid x \in B]$.  Then there exists a $K$-GDD of type $$\omega(\Pi):=\left[\sum_{x \in X_1} \omega(x),\dots,\sum_{x \in X_u} \omega(x)\right].$$
Furthermore, if the master GDD has strong dimension $d$, then the resultant GDD has strong dimension $\ge d$.
\end{cons}  

\proof
The construction proceeds by replacing each point $x \in X$ by a new set of $x_1,\dots,x_{\omega(x)}$ of $\omega(x)$ points, maintaining the group partition.  So the type becomes $\omega(\Pi)$.  Every block of the master, say $B \in \cB$, is replaced by a copy of the GDD of type $\omega(B)$ as defined.  In the resultant, if two points $x_i,y_j$ from different groups are given, their `projections' $x,y$ belong to different groups, and therefore a unique block in the master.  This block was replaced by a unique ingredient GDD.  It follows that $x_i,y_j$ appear together in a unique block in this ingredient, and therefore in the resultant.

For the claim on dimension, suppose a set of $d$ points is given in the resultant.  They arose from at most $d$ points, say $Y \subset X$ in the master GDD.  By assumption, $Y$ is contained in a proper strong subspace $X' \subset X$.  It is clear that $X'$ lifts to a strong proper subspace 
$\{x_i: x \in X', i=1,\dots,\omega(x)\}$ in the resultant, since two points from different groups in $X'$ lie on a block of some ingredient GDD placed on $X'$.
\qed

\Section{Proof of the main theorem}

Starting from the points and lines of AG$_d(q)$, let's apply Construction~\ref{wfc} with a large uniform weighting $\omega(x) \equiv n$, replacing blocks with TD$(q,n)$ for $n \ge n_0(q)$.  These ingredients exist by virtue of Theorem~\ref{asymptotic-td}.

\begin{prop}
For any positive integer $d$ and any prime power $q$, there exists a $q$-GDD of strong dimension $\ge d$ and type $n^{(q^d)} := [\underbrace{n,\dots,n}_{q^d}]$ for all sufficiently large integers $n$.
\end{prop}

Next, apply Construction~\ref{bub} to break up blocks of size $q$ by replacing them with PBD$(q,\{r\})$ when possible.  We can then truncate the last group via Construction~\ref{trunc}, dropping some block sizes by one.

\begin{prop}
\label{after-trunc}
For any positive integers $d$ and $r$ with $r \ge 3$, there exists an $\{r-1,r\}$-GDD of strong dimension $\ge d$ and type
$n^{(q^d-1)} x^1$
for all large prime powers $q \equiv 1 \pmod{r(r-1)}$, all sufficiently large integers $n$, and for any positive integer $x \le n$. 
\end{prop}

Now, take $r$ such that both $r$ and $r-1$ are replication numbers for PBDs with block size $k$.  That is, take $r$ such that $k$-GDD of type $(k-1)^r$ and $(k-1)^{r-1}$ both exist.  (It suffices to take $r \equiv 1 \pmod{k}$ and large, by Corollary~\ref{repl}, and delete a point.)    Apply Construction~\ref{wfc} once again, using weights $\omega(x) \equiv k-1$ and these $k$-GDDs.

\begin{prop}
For any positive integers $d$ and $k$ with $k \ge 2$, there exists a $k$-GDD of dimension $\ge d$ and type
$[n(k-1)]^{(q^d-1)} [x(k-1)]^1$
for large prime powers $q \equiv 1 \pmod{k(k-1)}$, all sufficiently large integers $n$, and for any positive integer $x \le n$. 
\end{prop}

\emph{Remark}.  The conditions on $q$ remain $q \equiv 1 \pmod{r(r-1)}$ and large; the preceding results hold for an infinite sequence of $q \equiv 1 \pmod{k(k-1)}$ since $k(k-1)|r(r-1)$.

Finally, invoke Lemma~\ref{add-pt}.  That is, add a point and fill groups with PBDs having block size $k$ and replication numbers $n,x$, which exist again by Corollary~\ref{repl} for admissible $n,x \ge r_0(k)$.  It is actually enough for our purposes to assume $k \mid n$.

\begin{prop}
\label{final-cons}
For any positive integers $d$ and $k$ with $k \ge 2$, there exists a PBD of blocksize $k$, dimension $\ge d$, and replication
number $n(q^d-1) + x$
for infinitely many prime powers $q \equiv 1 \pmod{k(k-1)}$, all sufficiently large integers $n$ with $k \mid n$, and for any integer $x$ with $r_0(k) \le x \le n$ and $x(x-1) \equiv 0 \pmod{k}$. 
\end{prop}

It remains to observe that these `intervals' of constructible replication numbers overlap for large $n$.  This is facilitated by the following easy observation.

\begin{lemma}
\label{overlap}
Given positive integers $A,c$, every sufficiently large integer $y$ with can be represented as
$y = n A + x$ for some integers $n$ and $x$, $c \le x \le n$.
\end{lemma}

\proof
Suppose $y \ge A(A+c+1)+c$ and apply the division algorithm to $y-c$ and $A$.  We have $y-c=nA+ m$, $0 \le m < A$.  Put $x=m+c$.
We have $A(A+c+1)+c \le y < (n+1)A+c$, which implies $A+c<n$.  It follows that $x$ lies in the required interval. 
\qed

We can now give a proof of (the weak version of) the main theorem.

{\sl Proof of Theorem~{\rm \ref{main-k}}}.
It suffices to prove that PBDs of blocksize $k$ and dimension $\ge d$ exist with all sufficiently large replication numbers $y$ satisfying $y(y-1) \equiv 0 \pmod{k}$.
Starting from the given $k$, let's choose $r$ and then $q$ as above.  Apply Lemma~\ref{overlap} with $A=k(q^d-1)$ and $c=r_0(k)$ to write $y=n(q^d-1)+x$ for some integer $n$ divisible by $k$, and where $r_0(k) \le x \le n$.  We may further assume that $y$ is sufficiently large so that $n \ge n_0(q)$.
Since $k \mid n$ and $x \equiv y \pmod{k(k-1)}$, the hypotheses of Proposition~\ref{final-cons} are satisfied.
This produces a PBD with the desired replication number $y$.
\qed

We turn our attention now to the full version, Theorem~\ref{main}.  For starters, it should be remarked that some (perhaps enough, after some work) of this follows as a corollary of Theorem~\ref{main-k}, appealing to Construction~\ref{bub}.  But it is easy enough to simply strengthen certain steps in the above proof.

Let $\alpha:=\alpha(K)$, $\beta:=\beta(K)$, and observe that $\alpha \mid \beta$.  Put $\gamma := \beta/\alpha$.
The necessary and asymptotically sufficient conditions for PBD$(v,K)$ can be rewritten as
$v=\alpha y + 1$, where $y$ is an integer satisfying $y(\alpha y+1) \equiv 0 \pmod{\gamma}$.  Note that this extends Corollary~\ref{repl}, since $\alpha y+1 \equiv 1-y \pmod{\gamma}$ in the case $\alpha=k-1$, $\gamma=k$.
As before, it suffices to realize PBDs of dimension $\ge d$ for all large values of $y$ having this form.  Let's suppose Wilson's theorem delivers a lower bound of $y_0=y_0(K)$  for existence of PBD$(\alpha y+1,K)$.

We work from Proposition~\ref{after-trunc}. Choose $r$ now so that $K$-GDD of type $\alpha^r$ and $\alpha^{r-1}$ both exist, appealing to Theorem~\ref{asymptotic-gdd}.  (It is enough to take $r \equiv 1 \pmod{\gamma}$ and large.)  Apply Construction~\ref{wfc} with constant weight $\alpha$. Fill with PBD$(n\alpha+1,K)$ and PBD$(x \alpha+1,K)$, which exist for sufficiently large $n,x$ in appropriate congruence classes.  For convenience, we can assume $\gamma \mid n$.  Here is the extension of Proposition~\ref{final-cons}. 

\begin{prop}
\label{final-cons-K}
For any positive integer $d$ and subset $K \subseteq \Z_{\ge 2}$, there exists a PBD$(v,K)$ of dimension $\ge d$ 
with $v=[n(q^d-1) + x]\alpha+1$ 
for infinitely many prime powers $q \equiv 1 \pmod{\beta}$, all sufficiently large integers $n$ with $\gamma \mid n$, and for any integer $x$ with $y_0 \le x \le n$ and $x(\alpha x+1) \equiv 0 \pmod{\gamma}$. 
\end{prop}

The strong version of our main result now easily follows.

{\sl Proof of Theorem~{\rm \ref{main}}}.
Let $y$ be large and admissible for PBD$(\alpha y +1,K)$.  Choose $r,q,n$ (in that order) according to the above requirements.  The choice of $n$ requires $n \ge n_0(q)$ but also comes from  Lemma~\ref{overlap} with $A=\gamma(q^d-1)$, and $c=y_0(K)$ so that $y=n(q^d-1)+x$ for $y_0 \le x \le n(q-1)$.  We have $\gamma \mid n$ and $x \equiv y \pmod{\beta}$.  It follows that $n,x$ are admissible for Proposition~\ref{final-cons-K},  yielding the desired PBD.
\qed

\Section{Discussion}

With similar (perhaps slightly more technical) constructions, it is often possible to compute an effective bound on $v$ for certain $d,K$.  For example, consider the case $d=3$, $K=\{3,4,5\}$.  This $K$ is interesting because $\alpha(K)=1$, $\beta(K)=2$ and in fact any positive integer $v \neq 2,6,8$ admits a PBD$(v,\{3,4,5\})$; see \cite{handbook} for instance.   In a recent thesis \cite{Niezen}, Niezen has shown that PBD$(v,\{3,4,5\})$ of dimension 3 exist for all $v \ge 48$, treating also many smaller values of $v$.

This brings up an interesting side note.  We have stated our main theorem with only a lower bound on dimension.  In concrete cases, such as $K=\{3,4,5\}$, it is possible to adapt an argument of Teirlinck in Section 2 of \cite{T} to `break' a space and reduce its dimension to exactly a desired $d$.  However, the argument relies on existence of (possibly small) spaces of dimension exactly two; while this is usually easy when such spaces exist, we are far from a complete existence theory for PBDs having general $K$.

Even more challenging is an existence theory for $t$-\emph{wise balanced designs}, where every $t$-element subset of points is contained in a unique block.  Our main theorem is not likely to prove useful for such objects, yet it does concern something weaker (but related).  In a PBD$(v,K)$ of dimension $t$, every $t$-subset of points is contained in a proper subdesign (which may be a single block).  These can be viewed as $t$-wise `covering' designs, but we must regard certain subdesigns as blocks.

Finally, we would like to mention a neat application of dimension and generated subspaces.  Consider again the case $d=3$ and $K=\{3,4,5\}$.  Suppose we slightly strengthen the dimension 3 requirement by (universally) bounding the three-point generated subspaces as $v$ grows.  (It should not be difficult to obtain a general result along these lines.)  For example, it was shown in \cite{DL} that, for all integers $v$, there exist PBD$(v,\{3,4,5\})$ such that any three points generate a subspace of size $< 1000$.  (This bound is far from best possible.)  An interesting consequence is that one can construct, using these linear spaces, one-factorizations (i.e. $n$-edge-colourings) of the complete bipartite graphs $K_{n,n}$ which universally bound the longest bi-colored cycle.  That this quantity can be universally bounded with respect to $n$ seems to be a surprising result.


\begin{thebibliography}{99}

\bibitem{BB}
L.M.~Batten and A.~Beutelspacher, The theory of finite linear spaces: combinatorics of points and lines, \emph{Cambridge University Press}, 1993.

\bibitem{CES}
S.~Chowla, P.~Erd\"os, \and E.G.~Strauss, On the maximal number of pairwise 
orthogonal Latin squares of a given order, {\em Canad. J. Math.} 12
(1960), 204--208.

\bibitem{handbook}
C.J.~Colbourn \and J.H.~Dinitz, eds., {\em The CRC Handbook of Combinatorial 
Designs}, 2nd edition, CRC Press, Inc., 2006.

\bibitem{D}
A.~Delandtsheer, 
Dimensional linear spaces.
{\it Handbook of incidence geometry}, 193--294, 
North-Holland, Amsterdam, 1995.

\bibitem{DL}
P.J.~Dukes \and A.C.H.~Ling, Linear spaces with small generated subspaces. {\em J. Combin. Theory A}, 116 (2009), 485--493.

\bibitem{LW}
E.R.~Lamken and R.M.~Wilson, Decompositions of edge-colored complete graphs.
{\em J. Combin. Theory Ser. A} 89 (2000), 149--200.

\bibitem{Liu}
J.~Liu, Asymptotic existence theorems for frames and group divisible designs. \emph{J. Combin. Theory Ser. A} 114 (2007), 410--420.

\bibitem{Niezen}
J.~Niezen, Pairwise balanced designs of dimension three. M.Sc. thesis, University of Victoria, 2013.

\bibitem{T}
L.~Teirlinck, On Steiner spaces, 
{\em J. Comb. Theory A} 26 (1979), 103--114.

\bibitem{W}
R.M.~Wilson, An existence theory for pairwise balanced designs: II, The 
structure of PBD-closed sets and the existence conjectures,
{\em J. Comb. Theory, Ser. A} 13 (1972), 246--273.

\end{thebibliography}
\end{document}